\documentclass[11pt]{amsart}
\usepackage{amsmath}
\usepackage{amsthm}
\usepackage{amssymb}
\usepackage{amscd}
\usepackage{bbold}
\usepackage[pdftex]{graphicx}
\usepackage{enumerate}

\newtheorem{lem}{Lemma}[section]
\newtheorem{thm}[lem]{Theorem}

\theoremstyle{definition}
\newtheorem{defn}[lem]{Definition}

\newcommand{\RR}{{\mathbb R}}
\newcommand{\Rhat}{\widehat{\mathbb{R}}}

\newcommand{\ZZ}{{\mathbb Z}}

\newcommand{\norm}[1]{\left\Vert #1\right\Vert}         
\newcommand{\inner}[2]{\left\langle #1, #2\right\rangle}    

\title{Balayage and short time Fourier \\transform frames}

\author{Enrico Au-Yeung}
\address{Pacific Institute for the Mathematical Sciences\\
Vancouver, BC Canada V6T 1Z4}
\email{enricoauy@math.ubc.ca} 

\author{John J. Benedetto}
\address{Norbert Wiener Center, Department of Mathematics \\ University of Maryland, College Park, MD 20742 USA}
\email{jjb@math.umd.edu}

\begin{document}

\begin{abstract}
Using his formulation of the potential theoretic notion of balayage
and his deep results about this idea, Beurling gave sufficient 
conditions for Fourier frames in terms of balayage. The analysis
makes use of spectral synthesis, due to Wiener and Beurling, as well
as properties of strict
multiplicity, whose origins go back to Riemann. In this setting and with
this technology, we formulate and prove
non-uniform sampling formulas in the context of the short time Fourier transform 
(STFT).

\end{abstract}

\maketitle



%


\section{Introduction}

\subsection{Background and theme}
Frames provide a natural tool for dealing with signal reconstruction in
the presence of noise in the setting of overcomplete sets of atoms, and with the goals
of numerical stability and robust signal representation. Fourier frames were 
originally studied in the 
context of non-harmonic Fourier series by Duffin and Schaeffer \cite{DufSch1952}, 
with a history going back to Paley and Wiener \cite{PalWie1934} (1934) 
and farther, and with significant activity in the 1930s and 1940s, e.g., see 
\cite{levi1940}. Since \cite{DufSch1952}, there have been significant contributions
 by Beurling (unpublished 1959-1960 lectures), \cite{beur1989}, \cite{beur1966}, Beurling and Malliavin 
\cite{BeuMal1962}, \cite{BeuMal1967}, Kahane \cite{Kaha1970}, Landau \cite{land1967}, Jaffard \cite{jaff1991}, and Seip \cite{Seip1995}, \cite{Seip2002}.

\begin{defn}(Frame)  Let $H$ be a separable Hilbert space.  A sequence $\{x_{n}\}_{n \in \ZZ} \subseteq H$ is a \emph{frame} 
for $H$ if there are positive constants $A$ and $B$ such that
\[\forall \ f \in H, \quad A \lVert f \rVert^{2} \leq \sum_{n \in \ZZ} |\langle f,x_{n}\rangle|^{2} \leq B \lVert f \rVert^{2} .  \]
The constants $A$ and $B$ are lower and upper frame bounds, respectively.  \end{defn}


Our overall goal is to formulate a general theory of Fourier frames and non-uniform sampling formulas 
parametrized by the space $M_{b}(\mathbb{R}^d)$ of bounded Radon measures, see \cite{BenAuY2013}. 
This formulation provides a natural way to generalize non-uniform sampling  to the setting of 
short time Fourier transforms (STFTs) \cite{groc2001}, Gabor theory \cite{groc1991}, \cite{FeiSun2006}, \cite{LabWeiWil2004}, and pseudo-differential operators \cite{groc2001}, \cite{groc2006},
 The techniques are based on Beurling's methods from 1959-1960, 
\cite{beur1966}, \cite{beur1989},  which incorporate balayage, spectral synthesis, and strict multiplicity. 
In this short paper, we show how to achieve this goal for STFTs.

\subsection{Definitions}
We define the Fourier transform $\mathcal{F}(f)$ of $f \in L^2(\mathbb{R}^d)$ and its inverse Fourier transform $\mathcal{F}^{-1}(f)$ by
\[\mathcal{F}(f) (\gamma) = \widehat{f}(\gamma) = \int_{\mathbb{R}^d} f(x) e^{-2 \pi i x \cdot \gamma} \ dx, \]  and  \[ \mathcal{F}^{-1}(\widehat{f})(\gamma) = f(x) = \int_{\mathbb{\widehat{R}}^d} \widehat{f}(\gamma) e^{2 \pi i x \cdot \gamma} \ d\gamma. \]
$\Rhat^d$ denotes $\RR^d$ considered as the spectral domain. We write $F^\vee(x) = \int_{\Rhat^d}F(\gamma)e^{2\pi i x \cdot \gamma}\,d\gamma$. The notation ``$\int$"' designates integration over $\RR^d$ or $\Rhat^d$. When $f$ is a bounded continuous function, its Fourier transform is defined in the sense of distributions.  If $X \subseteq \RR^d$, where $X$ is closed, then $M_b(X)$ is the space of bounded Radon measures $\mu$ with the support of $\mu$ contained in $X$. $C_b(\RR^d)$ denotes the space of complex-valued bounded continuous functions on $\RR^d$.

\begin{defn}(Fourier frame) Let $E \subseteq \mathbb{R}^d$ be a sequence and let $\Lambda \subseteq \widehat{\mathbb{R}}^d$ be a compact set. Notationally, let $e_{x}(\gamma) = e^{2 \pi i x \cdot \gamma}$. The sequence $\mathcal{E}(E) = \{e_{-x}: x \in E \}$ is a \emph{Fourier frame} for $L^2(\Lambda)$ if there are positive constants $A$ and $B$ such that
\begin{align*}
& \forall \ F \in L^2(\Lambda), \\
& A \lVert F \rVert^{2}_{L^2(\Lambda)} \leq \sum_{x \in E} |\langle F,e_{-x}\rangle|^{2} \leq B \lVert F \rVert^{2}_{L^2(\Lambda)}.  
\end{align*}
Define the \emph{Paley-Wiener space}, $$PW_{\Lambda}= \{f \in L^2(\RR^d): \textrm{supp} (\widehat{f}) \subseteq \Lambda\}.$$ Clearly, $\mathcal{E}(E)$ is a Fourier frame for $L^2(\Lambda)$ if and only if the sequence,
\[\{(e_{-x} \ \mathbb{1}_{\Lambda})^\vee: x \in E \} \subseteq PW_{\Lambda}, \] 
is a frame for $PW_{\Lambda}$, in which case it is called a \emph{Fourier frame} for $PW_{\Lambda}$. Note that $\inner{F}{e_{-x}} = f(x)$ for $f \in PW_{\Lambda}$, where $\widehat{f} = F \in L^2(\Rhat^d)$ can be  considered an element of $L^2(\Lambda).$
\end{defn}

Beurling introduced the following definition in his 1959-1960 lectures.
\begin{defn}(Balayage)
Let $E \subseteq \RR^d$ and $\Lambda \subseteq \Rhat^d$ be closed sets. \emph{Balayage} is possible for $(E, \Lambda) \subseteq \RR^d \times \Rhat^d$ if
\begin{equation*}
\forall \mu \in M_b(\RR^d), \mbox{ }\exists \nu \in M_b(E)  \mbox{ such that } \widehat{\mu} = \widehat{\nu} \mbox{ on } \Lambda .
\end{equation*}
\end{defn}
Balayage originated in potential theory, where it was introduced by Christoffel (early 1870s) and by Poincar{\'e} (1890). Kahane
formulated balayage for the harmonic analysis of restriction algebras. The set, $\Lambda$, of group characters (in this case $\RR^d$)
is the analogue of the original role of $\Lambda$ in balayage as a set of potential theoretic kernels.

Let $ \mathcal{C}(\Lambda) = \{f \in C_b(\RR^d):  \textrm{supp} (\widehat{f}) \subseteq \Lambda\}.$

\begin{defn}(Spectral synthesis)
A closed set $\Lambda \subseteq \Rhat^d$ is a set of \emph{spectral synthesis (S-set)} if
\begin{align*}
& \forall f \in \mathcal{C}(\Lambda) \text{ and } \forall \mu \in M_b(\mathbb{R}^d), \\
 & \widehat{\mu} = 0 \text{ on } \Lambda \Rightarrow \int f \,d\mu = 0,
\end{align*}
see \cite{bene1975}.
\end{defn}
Closely related to spectral synthesis is the ideal structure of $L^1$, which can be thought of as the Nullstellensatz of harmonic analysis.  
As examples of sets of spectral synthesis,
 polyhedra are S-sets, and the middle-third Cantor set is an S-set which contains non-S-sets.  Laurent Schwartz (1947) showed that $S^2 \subseteq \mathbb{\widehat{R}}^3$ is not an S-set; and, more generally, Malliavin (1959) proved that every non-discrete
locally compact abelian group contains non-S sets. See \cite{bene1975} for a unified treatment of this material.
\begin{defn}(Strict multiplicity)
A closed set $\Gamma \subseteq \Rhat^d$ is a set of \emph{strict multiplicity} if
\begin{equation*}
\exists \mu \in M_b(\Gamma)\setminus\{0\} \mbox{ such that } \lim_{\norm{x} \to \infty} |\mu^\vee (x) | = 0.
\end{equation*}
\end{defn}

The notion of strict multiplicity was motivated by Riemann's study of sets of uniqueness for trigonometric series.
Menchov (1906) showed that there exists a closed set $\Gamma \subseteq \widehat{\mathbb{R}}/ \mathbb{Z}$ and $\mu \in M(\Gamma) \setminus \{0\}$, such that
$| \Gamma | = 0$ and $\mu^{\vee}(n) = O((\log |n|)^{-1/2}), |n| \rightarrow \infty$. There have been intricate refinements
of Menchov's result by Bary (1927), Littlewood (1936), Beurling, et al., see \cite{bene1975}.

The above concepts are used in the deep proof of the following theorem.

\begin{thm}\label{theorem:balayage1}
 Assume that $\Lambda$ is an S-set of strict multiplicity, and that balayage is possible for $(E, \Lambda).$  Let $\Lambda_\epsilon = \{ \gamma \in \Rhat^d: \text{dist}\,(\gamma,\Lambda) \leq \epsilon \}$. There is ${\epsilon}_0 > 0$ such that
if $0 < \epsilon < {\epsilon}_0$, then 
 balayage is possible for $(E, \Lambda_\epsilon)$.

\end{thm}

\begin{defn}
 A sequence $E \subseteq \RR^d$ is \emph{separated} if $$\exists \, r >0 \text{ such that } \inf \{\norm{x-y}: x, y \in E \text{ and } x \neq y \} \geq r.$$
\end{defn}
The following theorem, due to Beurling, gives a sufficient condition for Fourier frames in terms of balayage.
Its history and structure are analyzed in \cite{BenAuY2013} as part of a more
general program. Theorem \ref{theorem:balayage1} is used in its proof.
\begin{thm}\label{theorem:balayage2}
 Assume that $\Lambda \subseteq \Rhat^d$ is an S-set of strict multiplicity and that $E \subseteq \RR^d$ is a separated sequence. If balayage is possible for $(E,\Lambda)$, then $\mathcal{E}(E)$ is a Fourier frame for $L^2(\Lambda)$, i.e., $\{(e_{-x} \ \mathbb{1}_{\Lambda})^\vee: x \in E \}$ is a Fourier frame for $PW_{\Lambda}$. 
\end{thm}
See \cite{land1967}, \cite{BenWu1999}, \cite{BenWu2_1999} (SampTA 1999), and \cite{BenWu2000}.




\section{Short time Fourier transform (STFT) frame inequalities}

\begin{defn}Let $f, g \in L^2(\mathbb{R}^{d})$.
 The {\it short time Fourier transform} (STFT) of $f$ with respect to $g$ is the function $V_{g}f$ on $\mathbb{R}^{2d}$ defined as
\[ \quad V_{g}f(x, \omega) = \int f(t) \overline{g(t-x)} \ e^{- 2 \pi i t \cdot \omega } \ dt,\]
\end{defn}
\noindent
see \cite{groc2001}, \cite{groc2006} (chapter 8).

The STFT is uniformly continuous on $\mathbb{R}^{2d}$.  Further, for a fixed ``window'' $g \in L^{2}(\mathbb{R}^{d})$ with $\|g\|_{2} = 1$,
we can recover the original function $f \in L^{2}(\mathbb{R}^{d})$ from its STFT $V_{g}f$ by means of the vector-valued integral inversion formula,
\begin{equation*}\label{eqn:InversionSTFT}
f = \int \int  V_{g}f(x, \omega) \ e_{\omega}{\tau}_{x} g \ d\omega \ dx,
\end{equation*}
\noindent
where $({\tau}_{x}g)(t) = g(t - x).$

\begin{thm}\label{theorem:semidiscrete}
Let $E = \{x_n\}\subseteq \RR^d$ be a
separated sequence, that is symmetric about $0 \in \mathbb{R}^{d}$; and let $\Lambda \subseteq \mathbb{R}^{d}$ be an $S$-set of strict multiplicity that is compact, convex, and symmetric about $0 \in \Rhat^d$. Assume balayage is possible for $(E, \Lambda.)$
Further, let $g \in L^2(\mathbb{R}^{d}), \widehat{g} = G,$ have the property that  $\norm{g}_2 = 1$.  

a. We have that
$$\quad \exists \ A > 0, \quad \text{such that } \quad \forall f \in PW_\Lambda \backslash \{0\}, \widehat{f} = F,$$
\begin{align*}
A \| f \|_{2}^{2} & \leq   \sum_{x \in E} \int | V_{G} F(\omega, x)|^{2} \ d{\omega} \\
& = \sum_{x \in E} \int | V_{g}f(x, \omega)|^2 \ d\omega.
\end{align*}

b. Let $G_{0}{}(\lambda) = 2^{d/4} e^{- \pi \| \lambda \|^2}$ so that $\| G_0 \|_{2} = 1$; and assume $\|V_{G_0}G \|_{1} < \infty$.  We have that 
$$\quad \exists \ B > 0, \quad \text{such that } \quad \forall f \in PW_\Lambda \backslash \{0\}, \widehat{f} = F,$$
\begin{align*}
\sum_{x \in E} \int | V_{g}f(x, \omega)|^2 \ d\omega & = \sum_{x \in E} \int | V_{G} F(\omega, -x)|^{2} \ d{\omega} \\
& \leq B \| f \|_{2}^{2},
\end{align*}
where $B$ can be taken as $C \| V_{G_0}G \|_{1}$ and where
$$\quad  C = \sup_{y, \gamma} \ \left\{ \sum_{x \in E} \int \left|V_{G_0}G_{0}(\gamma + \omega, y + x) \right| \ d\omega  \right\}.$$
\end{thm}


The technique of using $G_0$ goes back to Feichtinger and Zimmermann
\cite{FeiZimm1998}
(Lemma 3.2.15) for a related type of problem,
see also \cite{FeiSun2006} (Lemma 3.2).

We next consider balayage being possible for $(E, \Lambda)$, where $E = \{(s_m, t_n)\} \subseteq \RR^{2d}$ and $\Lambda \subseteq \Rhat^{2d}$.  This allows us to express the STFT $V_{g}f$ of $f$ as
$$
V_{g}f(y, \omega) = \sum_{m} \sum_{n}  a_{mn}(y, \omega) h(s_m - y, t_n - \omega) V_{g}f(s_m, t_n), 
$$
where 
$$ 
\sum_{m} \sum_{n} | a_{mn}(y, \omega) | < \infty.
$$
The  following result and others like it, including Theorem \ref{theorem:semidiscrete}, can be formulated in terms of $(X,\mu)$ frames, 
\cite{Antoine2000}, \cite{GabardoHan2003}. \cite{ForRau2005}. 

\begin{thm}\label{thm_nonuniform_Gabor_frame}
Assume balayage is possible for $(E, \Lambda)$, where $E = \{(s_m, t_n)\} \subseteq \RR^{2d}$ is separated, and $\Lambda \subseteq \Rhat^{2d}$ is an S-set that is compact, convex, and symmetric about $0 \in  \Rhat^{2d}$.
Fix a window function $g \in L^{2}(\mathbb{R}^{d})$ such that $\| g \|_2 = 1.$ There are constants $A,B>0$, such that if $f \in L^2(\mathbb{R}^d)$ satisfies
the conditions,\\
$(1) \ V_{g}f \in L^{1}(\mathbb{R}^{2d})$ and\\
$(2) \ \mathcal{F}(V_{g}f)(\zeta_1, \zeta_2)$ has support $\subseteq \Lambda \subseteq \Rhat^{2d}$,

\noindent then
\begin{align*}
A \int |f(x)|^2 \ dx  &\leq  \sum_{m} \sum_{n} | V_{g}f(s_m, t_n) |^2 \\
& \leq B \int |f(x)|^2 \ dx. 
\end{align*}
\end{thm}
The hypothesis that $V_{g}f \in L^{1}(\mathbb{R}^{2d})$ means that $f$ belongs to the Feichtinger algebra $\mathcal{S}_0(\mathbb{R}^d)$.   It is the smallest Banach space that is invariant under translations and modulations. There are other equivalent characterizations of $\mathcal{S}_0(\mathbb{R}^d)$, see \cite{Fei1981}, \cite{FeiZimm1998}.    Fix a function $ \mathcal{S}_0(\mathbb{R}^d)$ and define the vector space $\mathcal{M}_1^1$ of all non-uniform Gabor expansions
\[f = \sum_{n=1}^{\infty} c_n \tau_{x_n} e_{\omega_n} g, \]
where $\{(x_n, \omega_n) \in \mathbb{R}^{2d}, n \in \mathbb{N}\}$ is an arbitrary countable set of numbers and $\sum_{n=1}^{\infty} |c_n| < \infty$. For this space, the norm is taken to be $\inf \sum_{n=1}^{\infty} |c_n|$, where the infimum is taken over all possible representations.   Then the vector space $\mathcal{M}_1^1$ coincides with $ \mathcal{S}_0(\mathbb{R}^d)$.  For functions in $\mathcal{S}_0(\mathbb{R}^d)$, Theorem \ref{thm_nonuniform_Gabor_frame} should be compared to the following theorem of 
Gr{\"o}chenig \cite{groc1991}, \cite{groc2001} (Chapter 12):
\begin{thm}
\label{thm:groch}
Given any $g \in  \mathcal{S}_0(\mathbb{R}^d)$. There is $r = r(g) > 0$ such that if
$E = \{(s_n, \sigma_n)\} \subseteq {\mathbb R}^d \times {\widehat{\mathbb R}}^d$ is a separated sequence with the property that
$$
   \bigcup_{n=1}^{\infty} \overline{B((s_n,{\sigma}_n),r(g))} ={\mathbb R}^d \times {\widehat{\mathbb R}}^d,
$$
then the frame operator, $S = S_{g,E},$ defined by
$$
   S_{g,E}\,f = {\sum}_{n=1}^{\infty}\langle f, {\tau}_{s_n}e_{\sigma_n}g\rangle\, {\tau}_{s_n}e_{\sigma_n}g,
$$
is invertible on $\mathcal{S}_0(\mathbb{R}^d)$.

Moreover, every $f \in  \mathcal{S}_0(\mathbb{R}^d)$ has a non-uniform Gabor expansion,
$$
f = {\sum}_{n=1}^{\infty} \langle f, \tau_{x_n} e_{\omega_n} g \rangle S_{g,E}^{-1}(\tau_{x_n} e_{\omega_n}g),
$$
where the series converges unconditionally in $ \mathcal{S}_0(\mathbb{R}^d)$. \\
($E$ depends on $g.$)

\end{thm}

A critical, thorough comparison of Theorems II.3 and II.4 is
given in \cite{BenAuY2013}.

\section*{Acknowledgment}
The first named author gratefully acknowledges the support of MURI-AFOSR Grant
FA9550-05-1-0443. The second named author  gratefully acknowledges the support of MURI-ARO
Grant W911NF-09-1-0383 and NGA Grant HM-1582-08-1-0009. Both authors also benefitted
from insightful observations by Professors Carlos Cabrelli, Matei Machedon, Ursula Molter, and
Kasso Okoudjou, as well as from Dr. Henry J.~Landau, the grand master of Fourier frames.
Finally, we would like to thank the referees for their insightful observations.



\bibliographystyle{plain}
\bibliography{2012JBbibSAMPTA2013}

\begin{thebibliography}{10}

\bibitem{Antoine2000}
Syed~Twareque Ali, Jean-Pierre Antoine, and Jean-Pierre Gazeau.
\newblock {\em Coherent states, wavelets and their generalizations}.
\newblock Graduate {T}ext in {C}ontemporary {P}hysics. Springer-Verlag, New
  York, 2000.

\bibitem{BenWu2_1999}
J.~J. Benedetto and H.~Wu.
\newblock A {B}eurling covering theorem and multidimensional irregular
  sampling.
\newblock In {\em SampTA}, Loen, 1999.

\bibitem{bene1975}
John~J. Benedetto.
\newblock {\em Spectral {S}ynthesis}.
\newblock Academic Press, Inc., New York-London, 1975.

\bibitem{BenAuY2013}
John~J. Benedetto and Enrico Au-Yeung.
\newblock Generalized {F}ourier frames in terms of balayage.
\newblock {\em forthcoming}.

\bibitem{BenWu1999}
John~J. Benedetto and Hui-Chuan Wu.
\newblock A multidimensional irregular sampling algorithm and applications.
\newblock {\em IEEE-ICASSP}, 1999.

\bibitem{BenWu2000}
John~J. Benedetto and Hui-Chuan Wu.
\newblock Non-uniform sampling and spiral {M}{R}{I} reconstruction.
\newblock {\em SPIE}, 2000.

\bibitem{BeuMal1962}
A.~Beurling and P.~Malliavin.
\newblock On {F}ourier transforms of measures with compact support.
\newblock {\em Acta Mathematica}, 107:291--309.

\bibitem{beur1966}
Arne Beurling.
\newblock Local harmonic analysis with some applications to differential
  operators.
\newblock {\em Some Recent Advances in the Basic Sciences, Vol. 1 (Proc. Annual
  Sci. Conf., Belfer Grad. School Sci., Yeshiva Univ., New York, 1962--1964)},
  pages 109--125, 1966.

\bibitem{beur1989}
Arne Beurling.
\newblock {\em The {C}ollected {W}orks of {A}rne {B}eurling. {V}ol. 2.
  {H}armonic {A}nalysis.}
\newblock Birkh{\"a}user, Boston, 1989.

\bibitem{BeuMal1967}
Arne Beurling and Paul Malliavin.
\newblock On the closure of characters and the zeros of entire functions.
\newblock {\em Acta Mathematica}, 118:79--93.

\bibitem{DufSch1952}
Richard~J. Duffin and A.~C. Schaeffer.
\newblock A class of nonharmonic {F}ourier series.
\newblock {\em Trans. Amer. Math. Soc.}, 72:341--366, 1952.

\bibitem{Fei1981}
Hans~G. Feichtinger.
\newblock On a new {S}egal algebra.
\newblock {\em Monatsh. Math.}, 92:269--289, 1981.

\bibitem{FeiSun2006}
Hans~G. Feichtinger and Wenchang Sun.
\newblock Stability of {G}abor frames with arbitrary sampling points.
\newblock {\em Acta Mathematica Hungarica}, 113:187--212, 2006.

\bibitem{FeiZimm1998}
Hans~G. Feichtinger and Georg Zimmermann.
\newblock A {B}anach space of test functions for {G}abor analysis.
\newblock In {\em Gabor analysis and algorithms,}, Appl. Numer. Harmon. Anal.,
  pages 123--170. Birkh{\"a}user, Boston, 1998.

\bibitem{ForRau2005}
Massimo Fornasier and Holger Rauhut.
\newblock Continuous frames, function spaces, and the discretization problem.
\newblock {\em J. Fourier Anal. Appl.}, 11(3):245--287, 2005.

\bibitem{GabardoHan2003}
Jean-Pierre Gabardo and Deguang Han.
\newblock Frames associated with measurable spaces. frames.
\newblock {\em Adv. Comput. Math.}, 18(2-4):127--147, 2003.

\bibitem{groc1991}
Karlheinz Gr{\"o}chenig.
\newblock Describing functions: atomic decompositions versus frames.
\newblock {\em Monatsh. Math.}, 112:1--42, 1991.

\bibitem{groc2001}
Karlheinz Gr{\"o}chenig.
\newblock {\em Foundations of {T}ime-{F}requency {A}nalysis}.
\newblock Applied and Numerical Harmonic Analysis. Birkh\"auser Boston Inc.,
  Boston, MA, 2001.

\bibitem{groc2006}
Karlheinz Gr{\"o}chenig.
\newblock A pedestrian's approach to pseudodifferential operators.
\newblock In {\em Harmonic {A}nalysis and {A}pplications}, Appl. Numer. Harmon.
  Anal., pages 139--169. Birk\"{a}user, Boston, 2006.

\bibitem{jaff1991}
St{\'e}phane Jaffard.
\newblock A density criterion for frames of complex exponentials.
\newblock {\em Michigan Math. J.}, 38:339--348, 1991.

\bibitem{Kaha1970}
J.-P. Kahane.
\newblock Sur certaines classes de s\'eries de {F}ourier absolument
  convergentes.
\newblock {\em J. Math. Pures Appl. (9)}, 35:249--259, 1956.

\bibitem{LabWeiWil2004}
D.~Labate, G.~Weiss, and E.~Wilson.
\newblock An {A}pproach to the {S}tudy of {W}ave {P}acket {S}ystems.
\newblock {\em Contemporary Mathematics, Wavelets, Frames, and Operator
  Theory}, 345:215--235, 2004.

\bibitem{land1967}
Henry~J. Landau.
\newblock Necessary density conditions for sampling and interpolation of
  certain entire functions.
\newblock {\em Acta Mathematica}, 117:37--52, 1967.

\bibitem{levi1940}
Norman Levinson.
\newblock {\em Gap and {D}ensity {T}heorems}, volume XXVI of {\em Amer. Math.
  Society Colloquium Publications}.
\newblock American Mathematical Society, Providence, RI, 1940.

\bibitem{Seip2002}
Joaquim Ortega-Cerd\`{a} and K.~Seip.
\newblock Fourier frames.
\newblock {\em Ann. of Math.}, 155(3):789--806, 2002.

\bibitem{PalWie1934}
Raymond~E.~A.~C. Paley and Norbert Wiener.
\newblock {\em Fourier {T}ransforms in the {C}omplex {D}omain}, volume XIX of
  {\em Amer. Math. Society Colloquium Publications}.
\newblock American Mathematical Society, Providence, RI, 1934.

\bibitem{Seip1995}
K.~Seip.
\newblock On the connection between exponential bases and certain related
  sequences in {L}$^2(-\pi,\pi)$.
\newblock {\em J. Funct. Anal.}, 130:131--160, 1995.

\end{thebibliography}
%
%
%

\end{document}